\documentclass[reqno]{amsart}
\usepackage{amsmath, amssymb, amsthm, geometry, enumerate, graphicx, amsfonts, hyperref, mathrsfs}
\hypersetup{colorlinks=true,linkcolor=red, anchorcolor=blue, citecolor=blue, urlcolor=red, filecolor=magenta, pdftoolbar=true}
\theoremstyle{plain}
\newtheorem{thm}{\bf Theorem}[section]
\newtheorem{prop}[thm]{\bf Proposition}
\newtheorem{cor}[thm]{\bf Corollary}

\theoremstyle{remark}
\newtheorem{defn}[thm]{\bf Def{}inition}
\newtheorem{rem}[thm]{\bf Remark}
\newtheorem{exa}[thm]{\bf Example}
\numberwithin{equation}{section}
\usepackage{enumerate}
\usepackage{graphicx}
\allowdisplaybreaks
\begin{document}
\baselineskip08pt
\title [Super Weaving Frames]{Vector-Valued (Super) Weaving Frames}

\author[Deepshikha]{Deepshikha \ }
\address{Deepshikha, Department of Mathematics,
University of Delhi, Delhi-110007, India}
\email{dpmmehra@gmail.com}

\author[Lalit K. Vashisht]{ Lalit  Kumar  Vashisht}
\address{Lalit  K. Vashisht, Department of Mathematics,
University of Delhi, Delhi-110007, India}
\email{lalitkvashisht@gmail.com}

\begin{abstract}
 Two frames $\{\phi_{i}\}_{i \in I}$  and $\{\psi_{i}\}_{i \in I}$ for a separable  Hilbert space $H$ are  woven if there are positive constants $A \leq B$ such that for every subset $\sigma \subset I$, the family $\{\phi_{i}\}_{i \in \sigma} \cup \{\psi_{i}\}_{i \in \sigma^{c}}$  is a frame for $H$ with frame bounds $A,  B$.  Bemrose et al. introduced weaving frames in separable Hilbert spaces  and observed that weaving frames has potential applications in signal processing. Motivated by this, and the  recent work  of  Balan in the direction of application of vector-valued frames (or superframes) in signal processing, we study vector-valued weaving frames. In this paper, first we give some fundamental properties of vector-valued weaving frames. It is shown that if a family of   vector-valued frames is woven, then the corresponding family of frames for atomic spaces is woven, but the  converse is not true. We present a technique for the construction of vector-valued  woven frames from  given  woven frames for atomic  spaces . Necessary and sufficient conditions for vector-valued weaving Riesz sequences are given.   Several numerical examples are given to illustrate the results.
\end{abstract}


\keywords{ Frames,   vector-valued frames, weaving frames,  Riesz basis.\\
The research of Deepshikha  is supported by CSIR, India vide File No.: 09/045(1352)/2014-EMR-I.  Lalit was partially   supported by R $\&$ D Doctoral Research Programme, University of Delhi,
Delhi-110007, India (Grant
No.: RC/2015/9677).}

\maketitle

\baselineskip12pt
\section{Introduction}
Duffin and Schaeffer \cite{DS} introduced frames for Hilbert spaces in the context of nonharmonic Fourier series. Today frames have applications in a wide range of areas in applied mathematics. The applications of frames in  signal processing  are
now well-known, for example see \cite{CK, OC1}.  Balan  \cite{B1} introduced the concept of a vector-valued  frame (or ``superframe'')   in the context of multiplexing and further  studied in \cite{B2}.  The vector-valued frame has significant applications in mobile communication,
satellite communication, and computer area network. Recently,  Bemrose,  Casazza,  Gr\"{o}chenig,   Lammers  and  Lynch in \cite{BCGLL} proposed weaving frames in a separable Hilbert space. The concept of weaving frames is motivated by  a problem regarding distributed signal processing where redundant building blocks (frames) plays an important role. For example, in wireless sensor networks where frames may be subjected to distributed processing under different frames. Motivated by the concept of weaving frames and  superframes  and their application in Gabor and wavelet analysis, in this paper,  we study vector-valued (super) weaving frames.  Weaving
frames has potential applications in wireless sensor networks that require distributed
processing under different frames, as well as preprocessing of signals using Gabor frames. Notable contribution in the paper is a new technique for the construction of vector-valued  weaving frames from frames of atomic spaces. Some  necessary and sufficient conditions for vector-valued weaving Riesz sequences are given. Finally, a result for vector-valued weaving Riesz sequences in terms of operators on atomic spaces has been obtained.

\subsection{Previous works on weaving frames}

For a positive integer  $m$, we write $[m] = \{1,2,\dots,m\}$. We start with the definition of weaving frames in separable Hilbert spaces.

\begin{defn}\cite{BCGLL}
A family of frames $\{\phi_{ij}\}_{i\in I}$ for $j \in [m]$ for a Hilbert space $H$ is said to
be \emph{woven} if there are universal constants $A$ and $B$, so that for every partition $\{\sigma_j\}_{j \in [m]}$ of $I$, the family $\{\phi_{ij}\}_{i\in\sigma_j, j \in [m]}$ is a frame for $H$ with lower and upper frame bounds $A$
and $B$, respectively.
\end{defn}
\begin{defn}\cite{BCGLL}\label{def2}
A family of frames $\{\phi_{ij}\}_{i\in\mathbb{N}, j \in [m]}$ for a Hilbert space $H$ is said to be  \emph{weakly woven} if for every partition $\{\sigma_j\}_{j \in [m]}$ of $\mathbb{N}$, the family $\{\phi_{ij}\}_{i\in\sigma_j, j \in [m]}$ is a frame for $H$.
\end{defn}
\begin{rem}\label{rem2.5}
It is proved in \cite{BCGLL} that this weaker form of weaving (given in Definition \ref{def2})
 is equivalent to weaving
\end{rem}
Bemrose et al. proved in \cite{BCGLL} that a frame (which is not a Riesz basis) cannot be weave with
a Riesz basis.
\begin{thm}\cite{BCGLL}\label{th2.5}
Let $\{\phi_{j}\}_{j=1}^\infty$  be a Riesz basis and let $\{\psi_{j}\}_{j=1}^\infty$ be a frame for $H$. If  $\{\phi_{j}\}_{j=1}^\infty$ and  $\{\psi_{j}\}_{j=1}^\infty$ are woven, then $\{\psi_{j}\}_{j=1}^\infty$  must actually be a Riesz basis.
\end{thm}
Bemrose et al. in \cite{BCGLL}  classified  when Riesz bases and Riesz basic sequences can be woven
and provide a characterization in terms of distances between subspaces. Furthermore, they proved that if two Riesz bases are woven, then every weaving is in fact a Riesz basis, and not just a frame. A geometric characterization of woven Riesz bases in terms of distance between subspaces of a Hilbert space $H$ is given in \cite{BCGLL}. Casazza and Lynch in \cite{CL} reviewed fundamental properties of weaving frames.  They considered  a relation of frames to projections  and gave
a better understanding of what it really means for two frames to be woven.  Finally, they discussed  a weaving equivalent
of an unconditional basis. Generalized weaving frames in separable Hilbert spaces were  studied in \cite{DVV, DVash, VGDD}. Some necessary and sufficient conditions for weaving fusion frames may found in \cite{DV17}. In \cite{VD, VD2},  Vashisht and Deepshikha studied weaving frames with respect to measure spaces.   D\"{o}rfler and  Faulhuber  \cite{MoMF} studied weaving Gabor frames in $L^2(\mathbb{R})$.  They also discussed a family of  localization operators related to weaving Gabor frames.  A sufficient criteria for
a family of multi-window Gabor frames to be woven may found in \cite{MoMF}. Casazza, Freeman and Lynch \cite{CFL} extended the concept of weaving Hilbert space frames to the Banach space setting. They introduced and studied  \emph{weaving Schauder frames} in Banach spaces.
\subsection{Overview and Main Results}
The  paper is organized as follows:  Section \ref{sec2} contains the basic definitions and results about frames, weaving frames and vector-valued frames.  In Section \ref{sec3}, we introduce vector-valued weaving frames. Some properties of vector-valued frames are given.  It is shown that  frames associated with atomic spaces are woven, provided the given family of vector-valued frames is woven, see Proposition \ref{prop3.5}. The  converse is not true. Theorem \ref{thm3.9m} gives a technique for the construction of vector-valued  woven frames from  woven frames for  atomic spaces.  In Section \ref{sec4}, we give necessary and  sufficient conditions for vector-valued weaving Riesz sequences. Proposition \ref{prop4.4} gives a sufficient condition for vector-valued weaving Riesz bases.
It  is shown in Theorem \ref{th4.7} that if a  family of frames associated with one of the atomic spaces is a woven Riesz sequence, then the  family of vector-valued frames is also a  woven Riesz sequence. Finally, sufficient conditions for vector-valued weaving Riesz sequences in terms of bounded operators on atomic spaces  are  given in Theorem \ref{th4.9}. Several examples and counter-example are given to illustrate our results.
\section{Preliminaries}\label{sec2}
In this section, we  review the concepts of frames, vector-valued frames and weaving frames.  We begin with some notations: The set of all positive integers is denoted by  $\mathbb{N}$ or $\mathbb{Z}^+$. By $\mathbb{Z}$ and $\mathbb{Z}^-$ we denote the set of all integers and negative integers, respectively and $\mathbb{Z}^* = \mathbb{Z} \setminus \{0\} $.  For a set $E$, $\chi_E$ denote the characteristic function of $E$. By  $\mathbb{K}$ we denote field of real or complex numbers. For an indexing set $I$,  as usual $\ell^2(I)$ is the space  of all square summable complex-valued sequences  indexed by the set $I$. The space of bounded linear operators from a normed space $X$ into a normed space $Y$ is denoted by $\mathcal{B}(X, Y)$. If $X = Y$, then we write $\mathcal{B}(X, Y) = \mathcal{B}(X)$.

\subsection{Hilbert space frames:}
A countable sequence $\{f_k\}_{k \in  I}$ in a separable Hilbert space $H$ is
called a \emph{frame } (or \emph{Hilbert frame}) for $H$ if there exist positive numbers $ \alpha_o \leq \beta_o < \infty$ such that
\begin{align}\label{1.1}
\alpha_o \|f\|^2\leq  \sum_{k \in I} |\langle f, f_k\rangle|^2 \leq \beta_o \|f\|^2
\ \text{for  all} \  f \in H.
\end{align}
The numbers $\alpha_o$ and $\beta_o$ are called \emph{lower} and \emph{upper frame bounds}, respectively.  If upper inequality in \eqref{1.1} is satisfied,   then we say that  $\{f_k\}_{k \in I}$ is  a \emph{Bessel sequence} (or \emph{Hilbert Bessel sequence}) with \emph{Bessel bound} $\beta_o$. The frame $\{f_k\}_{k \in I}$ is \emph{tight} if it is possible to choose $\alpha_o = \beta_o$.

Let $\{f_k\}_{k \in  I}$  be a frame (or a Bessel sequence) for $H$. The \emph{pre-frame operator} operator is a linear bounded operator $T: \ell^2(I) \rightarrow H$  given by
\begin{align*}
T\{c_k\}_{k \in I} =\sum_{k \in I}c_k f_k.
\end{align*}
The frame operator $S = TT^*: H \rightarrow H $ is
\begin{align*}
Sf =\sum_{k \in I}\langle f, f_k \rangle f_k .
\end{align*}
 The frame operator for a frame $\{f_k\}_{k \in  I}$ is bounded linear and invertible operator on $H$. This gives the \emph{reconstruction formula} for all $f \in H$,
\begin{align*}
&f = SS^{-1}f =\sum_{k \in I} \langle S^{-1}f, f_k \rangle f_k
   \ =\sum_{k \in I}\langle f,S^{-1}f_k\rangle f_k .
\end{align*}
The scalars $\{\langle S^{-1}f, f_k \rangle\}_{k \in I}$ are called \emph{frame coefficients} of the vector $f \in H$.  One may observe that the representation of
$f$ in the reconstruction formula need not be unique.

A \emph{Riesz basis} for $H$ is a family of the form $\{Ue_k\}_{k\in I}$,   where $\{e_k\}_{k\in I}$ is an orthonormal basis for $H$ and $U:H\rightarrow H$ is a bounded linear  bijective operator. The following  result for bases in Banach spaces can be found in \cite[p. 173]{CH}
\begin{thm}\cite{CH}\label{th2.1}
Let $\{x_n\}_{n \in \mathbb{N}}$ be a basis for a Banach space $X$. If $\{y_n\}_{n \in \mathbb{N}} \subset X$ and there exists a constant $\lambda \in [0, 1)$ such that
\begin{align*}
\Big\|\sum_{n=1}^{N} c_n(x_n-y_n) \Big\| \leq \lambda \ \Big\|\sum_{n=1}^{N} c_n x_n\Big\|, \ N \in \mathbb{N},  \ \ c_1, c_2,\dots,c_N \in \mathbb{K}.
\end{align*}
then $\{y_n\}_{n \in \mathbb{N}}$ is a basis for $X$ and $\{y_n\}_{n \in \mathbb{N}}$ is equivalent to $\{x_n\}_{n \in \mathbb{N}}$. That is, there exists a bounded, linear and bijective operator $T:X \rightarrow X$ such that $Tx_n = y_n$ for all $n \in \mathbb{N}$.
\end{thm}

The basic theory of frames can be found in  books of Casazza and Kutyniok \cite{CK} and  Christensen \cite{OC1}.

\subsection{Vector-valued Franes (Superframes):}
Let  $I$ be a countable index set and consider \break
 $(\mathcal{F}_1; \pi_1; I), \dots,(\mathcal{F}_L; \pi_L; I)$, $L$ indexed sets of vectors (not necessarily
from the same Hilbert space), where $\pi_k : I \rightarrow F_k$ is the corresponding indexing map. A collection of such countable
sets of vectors together with their corresponding indexing maps from a same index set is called a \emph{superset}. In short, we write $(F_1, . . . ., F_L)$ for a superset when an indexing by a same index set $I$ for each subset $F_k$ of vectors of some Hilbert space $\mathcal{H}_k$ (or a
bigger space $K_k)$ is fixed.\\
 We write
\begin{align*}
\mathcal{F} = \mathcal{F}_1 \oplus \cdots \oplus \mathcal{F}_L = \{f^1_i \oplus \cdots \oplus f_i^L: i \in I \}, \ f^k_i = \pi_k(i) \in \mathcal{F}_k.
\end{align*}
Recall that the space $\mathcal{H}_1\oplus . . . . \oplus\mathcal{H}_L$ is a Hilbert space with natural inner product:
\begin{align*}
\langle f_1 \oplus \cdots \oplus f_L, \  g_1 \oplus \cdots \oplus g_L \rangle = \sum_{i =1}^{L} \langle f_i, g_i \rangle_{\mathcal{H}_i},  \ \Big(f_i, g_i \in \mathcal{H}_i \ (1 \leq i  \leq L)\Big).
\end{align*}
\begin{defn}\cite{B1}
The superset $(F_1, \dots, F_L)$ is called a \emph{vector-valued frame} or \emph{super frame} if $\mathcal{F}$ is a frame for the space $\mathcal{H}_1\oplus \cdots\oplus\mathcal{H}_L$. That is,  if there exists finite positive numbers $A_0 \leq B_0$ such that  for every $h_k \in \mathcal{H}_k \ (1 \leq k\leq L)$, we have
\begin{align*}
A_0(\| h_1\|^2 + \cdots  + \|h_L \|^2) \leq \sum_{i \in I} \Big| \sum_{k=1}^{L} \langle h_k, f^k_i \rangle \Big|^2\leq B_0 (\| h_1\|^2 + \cdots  + \|h_L \|^2).
\end{align*}
\end{defn}
If  $I$ is the common indexing set (possibly countably infinite) and $F_i = \{f_{i k}\}_{k \in I}, \ i = 1, 2,\dots,L$, then a vector-valued family is represented by $\{F_1, \dots, F_L, I\}$. The spaces $\mathcal{H}_i$ are called \emph{atomic spaces} (or \emph{components}) of the space $\mathcal{H}_1\oplus \cdots \oplus\mathcal{H}_L$.

By using superframes, Balan  in \cite{B1} discussed a signal transmission problem (namely, multiplexing, or Multiple Access)
and  analyzed several encoding-decoding schemes (with an example).
In recent years, vector-valued wavelet and Gabor frames have attracted mathematician and engineering specialists, for example see  \cite{ D1, DJ, F, GL, GH,   LH, LZ} and reference therein.

\section{Vector-Valued Weaving Frames}\label{sec3}
We start this section by the definition of a vector-valued weaving frame.
\begin{defn}\label{defn3.1I}
A family of vector-valued frames  $\Big\{\{F_1^i,\dots,F_L^i,I\}: i\in[m]\Big\}$ for $\mathcal{H}_1\oplus \cdots \oplus\mathcal{H}_L$ is said to be \emph{woven} if there exist universal constants $A$ and $B$ such that for any partition $\{\sigma_i\}_{i\in[m]}$ of $I$, the family
\begin{align*}
\bigcup_{i\in[m]}\{F_1^i,\dots,F_L^i, \sigma_i\}
\end{align*}
is a vector-valued frame for $\mathcal{H}_1\oplus \cdots \oplus\mathcal{H}_L$ with  lower and upper frame bounds $A$ and $B$, respectively. Each family $\bigcup_{i\in[m]}\{F_1^i,\dots,F_L^i, \sigma_i\}$ is called a \emph{vector-valued weaving} (or simply \emph{weaving}).
\end{defn}
\begin{rem}
Throughout the paper we take $F_j \subset \mathcal{H}_j$ for $j \in [L]$. If $\{F_1,\dots,F_L, I\}$ is vector-valued frame (Bessel sequence) for $\mathcal{H}_1\oplus \cdots\oplus\mathcal{H}_L$, then  $F_j \ (j = 1, \dots,L)$ are frames (Bessel sequences) for atomic spaces $\mathcal{H}_j$.
\end{rem}
Regarding the existence of  vector-valued weaving frames, we have the following examples.
\begin{exa}\label{exm3.3m}
Let $L = m =2; \ I = \mathbb{N}$ and let  $\mathcal{H}_1=\mathcal{H}_2= \ell^2(\mathbb{N})$. \\
 Define $F_1=\{f_{1i}\}_{i\in I},G_1=\{g_{1i}\}_{i\in I} \subset \mathcal{H}_1$ as follows:
\begin{align*}
f_{1i} = &\begin{cases}
e_{j},  \ i = 4j-3 \ (j \in \mathbb{N})\\
 0, \ \text{otherwise},
\end{cases}
\quad \text{and} \quad g_{1i} = \begin{cases}
e_{j},  \ i = 4j-3,4j-1 \ (j \in \mathbb{N})\\
0, \ \text{otherwise}.
\end{cases}
\end{align*}
and $F_2=\{f_{2i}\}_{i\in I},G_2=\{g_{2i}\}_{i\in I} \subset \mathcal{H}_2$ as follows:
\begin{align*}
f_{2i} = &\begin{cases}
e_{j},  \ i = 4j \ (j \in \mathbb{N})\\
0 \ \text{otherwise},
\end{cases}
\quad \text{and} \quad g_{2i} = \begin{cases}
e_{j},  \ i = 4j-2,4j \ (j \in \mathbb{N})\\
0 \ \text{otherwise},
\end{cases}
\end{align*}
where $\{e_i\}_{i \in I}$ is the  canonical orthonormal basis of $\ell^2(\mathbb{N}$.\\

To show the vector-valued families $\{F_1,F_2,I\}$ and $\{G_1,G_2,I\}$ are  woven, let  $\sigma \subset I$ be any subset.\\
We compute
\begin{align*}
&\sum_{i\in \sigma}\Big|\Big\langle f\oplus g,f_{1i}\oplus f_{2i}\Big\rangle\Big|^2+\sum_{i\in \sigma^c}\Big|\Big\langle f\oplus g,g_{1i}\oplus g_{2i}\Big\rangle\Big|^2\\
&\geq\sum_{i\in \sigma\bigcap\left(\{4m-3\}_{m\in\mathbb{N}}\cup\{4m\}_{m\in\mathbb{N}}\right)}\Big|\Big\langle f\oplus g,f_{1i}\oplus f_{2i}\Big\rangle\Big|^2
+\sum_{i\in \sigma^c\
\bigcap\left(\{4m-3\}_{m\in\mathbb{N}}\cup\{4m\}_{m\in\mathbb{N}}\right)}\Big|\Big\langle f\oplus g,g_{1i}\oplus g_{2i}\Big\rangle\Big|^2\\
&=\Big|\Big\langle f\oplus g,e_1\oplus 0\Big\rangle\Big|^2+\Big|\Big\langle f\oplus g,0\oplus e_1\Big\rangle\Big|^2+\Big|\Big\langle f\oplus g,e_2\oplus 0\Big\rangle\Big|^2 +\Big|\Big\langle f\oplus g,0\oplus e_2\Big\rangle\Big|^2+ \cdots\\
&=\Big|\langle f,e_1\rangle\Big|^2+\big|\langle g,e_1\rangle\Big|^2+\Big|\langle f,e_2\rangle\Big|^2+\Big|\langle g,e_2\rangle\Big|^2+ \cdots\\
&=\|f\|^2+ \|g\|^2\\
&=\|f\oplus g\|^2 \ \text{for all} \ f\oplus g\in\mathcal{H}_1\oplus\mathcal{H}_2.
\end{align*}
This gives a lower universal bound $A=1$. Similarly, we can show that an upper universal bound is $B=2$.
Hence the vector-valued frames  $\{F_1,F_2,I\}$ and $\{G_1,G_2,I\}$ for $\mathcal{H}_1 \oplus \mathcal{H}_2$ are  woven.
\end{exa}

\begin{exa}\label{ex3.4}
Let $L,  m, I, \ \mathcal{H}_1$ and $\mathcal{H}_2$ be same as in Example \ref{exm3.3m}.

 Define $F_1=\{f_{1i}\}_{i\in I},G_1=\{g_{1i}\}_{i\in I} \subset \mathcal{H}_1$ as follows:
\begin{align*}
f_{1i} = &\begin{cases}
e_{j},  \ i = 6j \ (j \in \mathbb{N})\\
 0, \ \text{otherwise}.
\end{cases}
\quad \text{and} \quad g_{1i} = \begin{cases}
e_{j},  \ i = 6j-3,6j-1,6j \ (j \in \mathbb{N})\\
0, \ \text{otherwise}.
\end{cases}
\end{align*}
and $F_2=\{f_{2i}\}_{i\in I},G_2=\{g_{2i}\}_{i\in I} \subset \mathcal{H}_2$ as follows:
\begin{align*}
f_{2i} = &\begin{cases}
e_{j},  \ i = 6j-1 \ (j \in \mathbb{N})\\
0 \ \text{otherwise}.
\end{cases}
\quad \text{and} \quad g_{2i} = \begin{cases}
e_{j},  \ i = 6j-4,6j-1,6j \ (j \in \mathbb{N})\\
0 \ \text{otherwise}.
\end{cases}
\end{align*}
where $\{e_i\}_{i\in I}$ is the  canonical orthonormal basis of $\mathcal{H}_1$. Then,  $\{F_{1},F_{2},I\}$ and $\{G_{1},G_{2},I\}$ are vector-valued frames for $\mathcal{H}_1\oplus\mathcal{H}_2$ but not woven.

 Indeed, suppose $\{F_{1},F_{2},I\}$ and $\{G_{1},G_{2},I\}$ are  woven  with universal bounds $A$ and $B$. \\
  For $\sigma = I\setminus\{5,6\}$ and $\left(-1,0,0,0, \dots \right)\oplus\left(1,0,0,0, \dots \right)\in\mathcal{H}_1\oplus\mathcal{H}_2$, we compute
\begin{align*}
&\sum_{i\in \sigma}\Big|\Big\langle \left(-1,0,0,0, \dots \right),f_{1i}\Big\rangle+\Big\langle \left(1,0,0,0, \dots \right),f_{2i}\Big\rangle \Big|^2\\
&\quad \quad \quad \quad \quad +\sum_{i\in \sigma^c}\Big|\Big\langle \left(-1,0,0,0, \dots \right),g_{1i}\Big\rangle+\Big\langle \left(1,0,0,0, \dots\right),g_{2i}\Big\rangle \Big|^2\\
&=\Big|\Big\langle \left(-1,0,0,0, \dots \right),e_1\Big\rangle+\Big\langle \left(1,0,0,0, \dots \right),e_1\Big\rangle \Big|^2+\Big|\Big\langle \left(-1,0,0,0, \dots \right),e_1\Big\rangle+\Big\langle \left(1,0,0,0, \dots \right),e_1\Big\rangle \Big|^2\\
&=0\\
&<2A =A\Big(\Big\|(-1 ,0,0,0,\dots ) \oplus(1,0,0,0, \dots )\Big\|^2\Big), \ \text{a contradiction}.
\end{align*}
 Hence the vector-valued frames $\{F_{1},F_{2},I\}$ and $\{G_{1},G_{2},I\}$ for $\mathcal{H}_1 \oplus \mathcal{H}_2$  are  not  woven.
\end{exa}

 The following proposition  shows that  if a family of vector-valued frames is woven, then it is componentwise woven. That is, frames associated with atomic spaces are woven. The converse is not true, see Example \ref{ex3.7}.
\begin{prop}\label{prop3.5}
Suppose  the family $\Big\{\{F_{1}^i,\dots, F_{L}^i, I\}:i\in[m]\Big\}$ of  vector-valued frames for \break $\mathcal{H}_1\oplus \cdots \oplus \mathcal{H}_L$ is woven with universal bounds $A$ and $B$. Then, for each $j\in[L]$, the  frames $\Big\{F_{j}^i:i\in[m]\Big\}$ for $\mathcal{H}_j$ are  woven.
\end{prop}
\proof
Let $j\in[L]$ be arbitrary but fixed and let $\{\sigma_i\} _{i\in [m]}$ be any partition of $I$. Let us write $F_p^i=\{f_{pk}^i\}_{k\in I}$,  $p\in[L], \ i \in [m]$.

 Then,  for any $f\in \mathcal{H}_j$  $($note that $0\oplus \cdots \oplus f\oplus \cdots \oplus0\in \mathcal{H}_1\oplus \cdots \oplus \mathcal{H}_j\oplus \cdots \oplus\mathcal{H}_L)$, we compute
\begin{align*}
A\|f\|^2
&\leq\sum_{i\in[m]}\sum_{k\in\sigma_i}\Big|\langle0,f_{1k}^i\rangle+ \cdots +\langle f,f_{jk}^i\rangle+ \cdots + \langle0,f_{Lk}^i\rangle\Big|^2\\
&=\sum_{i\in[m]}\sum_{k\in\sigma_i}\Big|\langle f,f_{jk}^i\rangle\Big|^2\\
&\leq B\|f\|^2.
\end{align*}
 Hence the family of frames $\Big\{F_{j}^i:i\in[m]\Big\}$  for the atomic space $\mathcal{H}_j$ is  woven.
\endproof
We now demonstrate by a concrete example that converse of Proposition  \ref{prop3.5} is not true.

\begin{exa}\label{ex3.7}
Let $F_1,  G_1 \subset \mathcal{H}_1$  and $ F_2,  G_2 \subset \mathcal{H}_2$ be frames for  $\mathcal{H}_1$ and  $\mathcal{H}_2$, respectively given in Example 3.4. Then, for any subset $\sigma$ of $I$, the families  $\{f_{1i}\}_{i\in\sigma}\bigcup\{g_{1i}\}_{i\in\sigma^{c}}$ and $\{f_{2i}\}_{i\in\sigma}\bigcup\{g_{2i}\}_{i\in\sigma^{c}}$ are frames for $\mathcal{H}_1$ and $\mathcal{H}_2$,  respectively  with  frame  bounds $1$ and $3$. That is, the frames $F_{i}$ and $G_{i}$ for  $\mathcal{H}_i$ $(i=1,2)$ are woven. But as shown in Example \ref{ex3.4}, the vector-valued frames $\{F_1, F_2, I\}$ and $\{G_1, G_2, I\}$ are not woven.
\end{exa}

However, this is not the case for Bessel sequences.
\begin{prop}\label{pro3.5}
Suppose that $\{F_{1}^i,\dots, F_{L}^i, I\}$  are vector-valued families  for  $\mathcal{H}_1\oplus \cdots \oplus\mathcal{H}_L$, where $F_k^i=\{f_{kj}^i\}_{j\in I}$ $(k\in[L], i\in[m])$ such that $\Big\{\{f_{kj}^i\}_{j\in I}:i\in[m]\Big\}$ is a woven Bessel sequence, for all $k\in[L]$ . Then, the family  $\Big\{\{F_{1}^i, \dots, F_{L}^i, I\}:i\in[m]\Big\}$ is a woven vector valued Bessel sequence for $\mathcal{H}_1\oplus \cdots \oplus\mathcal{H}_L$.
\end{prop}
\proof
Let $\Big\{\{f_{kj}^i\}_{j\in I}:i\in[m]\Big\}$ be a woven Bessel sequence with universal bound $B_k$, for all $k\in[L]$. Then, for any partition $\{\sigma_i\} _{i\in [m]}$ of $I$ and for any $(h_1, \dots, h_L) \in \mathcal{H}_1\oplus \cdots \oplus\mathcal{H}_L$, we have
\begin{align*}
\sum_{i\in [m]} \sum_{j\in\sigma_i}\Big|\sum_{k=1}^L\langle h_k,f_{kj}^i\rangle\Big|^2 &\leq\sum_{i\in [m]} \sum_{j\in \sigma_i}\sum_{k=1}^L 2^{L-1}\Big|\langle h_k,f_{kj}^i\rangle\Big|^2\\
&=\sum_{k=1}^L 2^{L-1}\left(\sum_{i\in [m]} \sum_{j\in \sigma_i}\Big|\langle h_k,f_{kj}^i\rangle\Big|^2\right)\\
&\leq\sum_{k=1}^L 2^{L-1}B_k\|h_k\|^2\\
& \leq 2^{L-1}\max_{k\in[L]}\{B_k\}\Big(\|h_1\|^2+ \cdots +\|h_L\|^2\Big).
\end{align*}
This concludes the proof.
\endproof

\begin{cor}\label{cor3.8}
Suppose that $\{F_{1}^i, \dots, F_{L}^i, I\}$  are vector-valued families  for  $\mathcal{H}_1\oplus \cdots \oplus\mathcal{H}_L$, where $F_k^i=\{f_{kj}^i\}_{j\in I}$ is a Bessel sequence for all $k\in[L], i\in[m]$. Then,   $\Big\{\{F_{1}^i, \dots, F_{L}^i, I\}:i\in[m]\Big\}$ is a woven vector-valued Bessel sequence for $\mathcal{H}_1\oplus \cdots \oplus\mathcal{H}_L$.
\end{cor}

A natural question arises about the  construction of vector-valued woven frames  from  given  woven frames for atomic spaces. In this direction,  the following theorem provides  a technique for  the construction of a family of vector-valued woven frames from given woven frames for atomic spaces. We prove the result for $L = 2 $. This can be extended  to  any finite superset.
\begin{thm}\label{thm3.9m}
 Assume that $\Phi_j = \{\phi_i^j\}_{i\in I_j}$ and  $\Psi_j = \{\psi_i^j\}_{i\in I_j}$   are woven frames for  the Hilbert space $\mathcal{H}_j$ $(j = 1, 2)$. Then, there is a vector-valued family associated with $\Phi_j$ and $\Psi_j$ $(j = 1,2)$  which constitutes vector-valued  woven frames for  $\mathcal{H}_1\oplus\mathcal{H}_2$.
\end{thm}
\proof Since $I_1$ and $I_2$ are countably infinite sets,  there exist bijective maps $\Theta_1:\mathbb{Z}^{-}\rightarrow I_1$ and $\Theta_2:\mathbb{Z}^{+}\rightarrow I_2$.\\
Define $F_1=\{f_{1i}\}_{i\in \mathbb{Z}^{*}} \subset \mathcal{H}_1$ and $F_2=\{f_{2i}\}_{i\in \mathbb{Z}^{*}} \subset \mathcal{H}_2$ as follows:
 \begin{align*}
f_{1i} = &\begin{cases}
\phi_{\Theta_1(i)}^1,\ i\in \mathbb{Z}^{-}\\
0,  \ i\in \mathbb{Z}^{+},
\end{cases}
\quad \text{and} \quad f_{2i} = \begin{cases}
0,\ i\in \mathbb{Z}^{-}\\
\phi_{\Theta_2(i)}^2,  \ i\in \mathbb{Z}^{+}.
\end{cases}
\end{align*}
 and define
 $G_1=\{g_{1i}\}_{i\in \mathbb{Z}^{*}} \subset \mathcal{H}_1$ and $G_2=\{g_{2i}\}_{i\in \mathbb{Z}^{*}} \subset \mathcal{H}_2$ by
\begin{align*}
g_{1i} = &\begin{cases}
\psi_{\Theta_1(i)}^1,\ i\in \mathbb{Z}^{-}\\
0,  \ i\in \mathbb{Z}^{+}
\end{cases}
\quad \text{and} \quad
g_{2i} = \begin{cases}
0,\ i\in \mathbb{Z}^{-}\\
\psi_{\Theta_2(i)}^2,  \ i\in \mathbb{Z}^{+},
\end{cases}
\end{align*}

To show that the vector-valued family  $\{F_1,F_2,\mathbb{Z}^{*}\}$ and $\{G_1,G_2,\mathbb{Z}^{*}\}$  are woven, let $\sigma$ be any subset of $\mathbb{Z}^{*}$ and $f\oplus g\in \mathcal{H}_1\oplus \mathcal{H}_2$ be arbitrary.\\
We compute
\begin{align}\label{eq3.1}
&\sum_{i\in \sigma}\Big|\Big\langle f\oplus g,f_{1i}\oplus f_{2i}\Big\rangle\Big|^2+\sum_{i\in \sigma^c}\Big|\Big\langle f\oplus g,g_{1i}\oplus g_{2i}\Big\rangle\Big|^2 \notag\\
& =\sum_{i\in \sigma}\Big|\Big\langle f,f_{1i}\Big\rangle+\Big\langle g,f_{2i}\Big\rangle \Big|^2+\sum_{i\in \sigma^c}\Big|\Big\langle f,g_{1i}\Big\rangle+\Big\langle g,g_{2i}\Big\rangle \Big|^2\notag\\
&=\sum_{i\in \sigma\cap \mathbb{Z}^{-}}\Big|\Big\langle f,f_{1i}\Big\rangle+\Big\langle g,f_{2i}\Big\rangle \Big|^2+\sum_{i\in \sigma\cap \mathbb{Z}^{+}}\Big|\Big\langle f,f_{1i}\Big\rangle+\Big\langle g,f_{2i}\Big\rangle \Big|^2+\sum_{i\in \sigma^c\cap \mathbb{Z}^{-}}\Big|\Big\langle f,g_{1i}\Big\rangle+\Big\langle g,g_{2i}\Big\rangle \Big|^2\notag\\
&\quad \quad \quad \quad+\sum_{i\in \sigma^c\cap \mathbb{Z}^{+}}\Big|\Big\langle f,g_{1i}\Big\rangle+\Big\langle g,g_{2i}\Big\rangle \Big|^2\notag\\
&=\sum_{i\in \sigma\cap \mathbb{Z}^{-}}\Big|\Big\langle f,f_{1i}\Big\rangle+\Big\langle g,0\Big\rangle \Big|^2+\sum_{i\in \sigma\cap \mathbb{Z}^{+}}\Big|\Big\langle f,0\Big\rangle+\Big\langle g,f_{2i}\Big\rangle \Big|^2+\sum_{i\in \sigma^c\cap \mathbb{Z}^{-}}\Big|\Big\langle f,g_{1i}\Big\rangle+\Big\langle g,0\Big\rangle \Big|^2\notag \\
& \quad \quad \quad \quad +\sum_{i\in \sigma^c\cap \mathbb{Z}^{+}}\Big|\Big\langle f,0\Big\rangle+\Big\langle g,g_{2i}\Big\rangle \Big|^2\notag\\
&=\sum_{i\in \sigma\cap \mathbb{Z}^{-}}\Big|\Big\langle f,\phi_{\Theta^1(i)}^1\Big\rangle\Big|^2+\sum_{i\in \sigma\cap \mathbb{Z}^{+}}\Big|\Big\langle g,\phi_{\Theta^2(i)}^2\Big\rangle \Big|^2+\sum_{i\in \sigma^c\cap \mathbb{Z}^{-}}\Big|\Big\langle f,\psi_{\Theta^1(i)}^1\Big\rangle\Big|^2+\sum_{i\in \sigma^c\cap \mathbb{Z}^{+}}\Big|\Big\langle g,\psi_{\Theta^2(i)}^2\Big\rangle \Big|^2\notag\\
&=\left(\sum_{i\in \sigma\cap \mathbb{Z}^{-}}\Big|\Big\langle f,\phi_{\Theta^1(i)}^1\Big\rangle\Big|^2+\sum_{i\in \sigma^c\cap \mathbb{Z}^{-}}\Big|\Big\langle f,\psi_{\Theta^1(i)}^1\Big\rangle\Big|^2\right) +\left(\sum_{i\in \sigma\cap \mathbb{Z}^{+}}\Big|\Big\langle g,\phi_{\Theta^2(i)}^2\Big\rangle \Big|^2+\sum_{i\in \sigma^c\cap \mathbb{Z}^{+}}\Big|\Big\langle g,\psi_{\Theta^2(i)}^2\Big\rangle \Big|^2\right)\notag\\
&=\left(\sum_{i\in \Theta^1(\sigma\cap \mathbb{Z}^{-})}\Big|\Big\langle f,\phi_i^1\Big\rangle\Big|^2+\sum_{i\in \Theta^1(\sigma^c\cap \mathbb{Z}^{-})}\Big|\Big\langle f,\psi_i^1\Big\rangle\Big|^2\right)\notag\\
& \quad \quad \quad \quad +\left(\sum_{i\in \Theta^2(\sigma\cap \mathbb{Z}^{+})}\Big|\Big\langle g,\phi_i^2\Big\rangle \Big|^2+\sum_{i\in \Theta^2(\sigma^c\cap \mathbb{Z}^{+})}\Big|\Big\langle g,\psi_i^2\Big\rangle \Big|^2\right).
\end{align}
Let $\sigma_1=\Theta^1(\sigma\cap \mathbb{Z}^{-}),\sigma_2=\Theta^2(\sigma\cap \mathbb{Z}^{+})$. Then,  $\sigma_1$ is a subset of $I_1$ such that the compliment of $\sigma_1$ in $I_1$ is  $\Theta^1(\sigma^c\cap \mathbb{Z}^{-})$ and $\sigma_2$ is a subset of $I_2$ such that the  compliment of $\sigma_2$ in $I_2$ is  $\Theta^2(\sigma^c\cap \mathbb{Z}^{+})$. By  using the fact that $\Phi_i$ and  $\Psi_i$ are woven frames for the  Hilbert space $\mathcal{H}_i \ (i = 1, 2)$ with universal bounds $A_i, B_i$ (say)  and using \eqref{eq3.1}, we compute
 \begin{align*}
& \min\{A_1,A_2\}(\|f\|^2+\|g\|^2)\\
&\leq A_1\|f\|^2+A_2\|g\|^2\\
& \leq \left(\sum_{i\in\sigma_1}\Big|\Big\langle f,\phi_i^1\Big\rangle\Big|^2+\sum_{i\in I_1\setminus\sigma_1}\Big|\Big\langle f,\psi_i^1\Big\rangle\Big|^2\right)+\left(\sum_{i\in\sigma_2}\Big|\Big\langle g,\phi_i^2\Big\rangle \Big|^2+\sum_{i\in I_2\setminus\sigma_2}\Big|\Big\langle g,\psi_i^2\Big\rangle \Big|^2\right)\\
& = \sum_{i\in \sigma}\Big|\Big\langle f\oplus g,f_{1i}\oplus f_{2i}\Big\rangle\Big|^2+\sum_{i\in \sigma^c}\Big|\Big\langle f\oplus g,g_{1i}\oplus g_{2i}\Big\rangle\Big|^2\\
&\leq B_1\|f\|^2+ B_2\|g\|^2\\
&\leq \max\{B_1,B_2\}(\|f\|^2+\|g\|^2).
\end{align*}
Hence the vector-valued frames $\{F_1,F_2,\mathbb{Z}^{*}\}$ and $\{G_1,G_2,\mathbb{Z}^{*}\}$  are woven with universal frame bounds $ \min\{A_1,A_2\}$ and $ \max\{B_1,B_2\}$.
\endproof

\section{Vector-valued Weaving Riesz Bases}\label{sec4}

\begin{defn}
The family  $\{F_1, \dots, F_L, I\}$ $(F_i=\{f_{ik}\}_{k\in I}\subseteq\mathcal{H}_i \ (1 \leq i \leq L))$ is called a \emph{vector-valued Riesz basis} (or \emph{ super-Riesz basis}) for the space $\mathcal{H}_1\oplus \cdots \oplus\mathcal{H}_L$ if there exists a  bounded bijective operator $U:\mathcal{H}_1\oplus \cdots \oplus\mathcal{H}_L\rightarrow \mathcal{H}_1\oplus \cdots \oplus\mathcal{H}_L$ such that $U
(e_{1k}\oplus \cdots \oplus e_{Lk})=f_{1k}\oplus \cdots\oplus f_{Lk}$ $(k\in I)$,   where $\{e_{1k}\oplus \cdots \oplus e_{Lk}\}_{k\in I}$ is an orthonormal basis for $\mathcal{H}_1\oplus \cdots \oplus\mathcal{H}_L$.
\end{defn}
A characterization of vector-valued  Riesz basis for a Hilbert space $\mathcal{H}_1\oplus . . . . \oplus\mathcal{H}_L$ can be obtained by using   \cite[Theorem 3.6.6]{OC1}
\begin{thm}
For a sequence $\{F_1, \dots, F_L,I\}$ $(F_i=\{f_{ik}\}_{k\in I}\subseteq\mathcal{H}_i)$, the following conditions are equivalent:
\begin{enumerate}[$(i)$]
\item $\{F_1, \dots, F_L,I\}$ is a vector-valued Riesz basis for $\mathcal{H}_1\oplus \cdots \oplus\mathcal{H}_L$.
\item $\{F_1, \dots, F_L,I\}$ is complete in $\mathcal{H}_1\oplus \cdots \oplus\mathcal{H}_L$ and there exist finite constants $A, B>0$ such that for every finite scalar sequence $\{c_k\}$, one has
\begin{align}\label{eq3.2}
A\sum|c_k|^2\leq \Big\|\sum c_k(f_{1k}\oplus \cdots\oplus f_{Lk})\Big\|^2 \leq B\sum|c_k|^2.
\end{align}
\end{enumerate}
\end{thm}
The scalars  $A$ and $B$  are called lower and upper \emph{ Riesz bounds}, respectively.
 \begin{defn}
The family $\{F_1, \dots, F_L,I\}$ $(F_i=\{f_{ik}\}_{k\in I}\subseteq\mathcal{H}_i)$ is called a \emph{vector-valued Riesz sequence} (or \emph{super-Riesz sequence}) in the space $\mathcal{H}_1\oplus \cdots \oplus\mathcal{H}_L$ if (4.1) holds for all finite scalar sequences $\{c_k\}_{k\in I}$.
\end{defn}

The following proposition  provides  a sufficient condition for vector-valued  Riesz basis.
\begin{prop}\label{prop4.4}
Suppose $\{F_1, \dots,F_L, I\}$ $(\text{where} \ F_i=\{f_{ik}\}_{k\in I})$ is a vector-valued  frame for $\mathcal{H}_1\oplus \cdots \oplus\mathcal{H}_L$ such that for some $n\in[L]$, the sequence $F_n$ is a Riesz basis for $\mathcal{H}_n$.  Then, the family  $\{F_1, \dots,F_L, I\}$ is a vector-valued  Riesz basis for $\mathcal{H}_1\oplus \cdots \oplus\mathcal{H}_L$.
\end{prop}
\proof
Suppose that $F_n$ is a Riesz basis for $\mathcal{H}_n$  with lower and upper Riesz bounds $A_n$ and $B_n$, respectively. Then, for any finite scalar sequence $\{c_k\}_{k\in I}$, we  compute
\begin{align*}
\Big\|\sum_{k\in I}c_k\Big(f_{1k}\oplus \cdots \oplus f_{Lk}\Big)\Big\|^2&=\Big\|\sum_{k\in I}c_kf_{1k}\Big\|^2+ \cdots +\Big\|\sum_{k\in I}c_kf_{Lk}\Big\|^2\\
&\geq\Big\|\sum_{k\in I}c_kf_{nk}\Big\|^2\\
&\geq A_n\sum_{k\in I}|c_k|^2.
\end{align*}
Therefore, \eqref{eq3.2} is satisfied.  By hypothesis, the family $\{F_1, \dots, F_L, I\}$ is a vector-valued frame for $\mathcal{H}_1\oplus \cdots \oplus\mathcal{H}_L$. Thus, the family $\{F_1, \dots,F_L, I\}$ is complete in $\mathcal{H}_1\oplus \cdots\oplus\mathcal{H}_L$. Hence  $\{F_1, \dots, F_L, I\}$ is a vector-valued Riesz basis for $\mathcal{H}_1\oplus \cdots \oplus\mathcal{H}_L$.
\endproof
\begin{rem}
If a  family $\{F_1, \dots,F_L, I\}$ is a vector-valued Riesz basis for $\mathcal{H}_1\oplus \cdots\oplus\mathcal{H}_L$, then its component $F_k$ need not be a Riesz basis for the atomic  space $\mathcal{H}_k$.  Indeed, let  $L =2, I= \mathbb{N}$ and $\mathcal{H}_1=\mathcal{H}_2= L^2(I, \mu)$, where $\mu$ is the counting measure. Define $F_{1}=\{f_{1i}\}_{i \in I},F_{2}=\{f_{2i}\}_{i\in I}$ as follows :
\begin{align*}
f_{1i} = &\begin{cases}
 e_j,\ i=2j-1\ (j \in \mathbb{N})\\
0,  \ i = 2j \ (j \in \mathbb{N}).
\end{cases}
;\quad  \quad f_{2i} = \begin{cases}
 0,\ i=2j-1\ (j \in \mathbb{N})\\
e_j,  \ i = 2j \ (j \in \mathbb{N}).
\end{cases}
\end{align*}
where $\{e_i\}_{i\in I}$ is the  canonical orthonormal basis of $\mathcal{H}_1$.\\
Then, $\{F_1,F_2, I\}$ is a vector-valued Riesz basis for $\mathcal{H}_1\oplus\mathcal{H}_2$. But $F_1$ and $F_2$ are not Riesz bases for $\mathcal{H}_1$ and $\mathcal{H}_2$, respectively.
\end{rem}

\begin{defn}
A  vector-valued family   $\Big\{\{F_1^i, \dots, F_L^i,I\}: i\in[m]\Big\}$ $(F_j^i=\{f_{jk}^i\}_{k\in I})$ for \break $\mathcal{H}_1\oplus \cdots \oplus\mathcal{H}_L$ is said to be a \emph{woven Riesz sequence} if there exist universal constants $A$ and $B$ such that for any partition $\{\sigma_i\}_{i\in[m]}$ of $I$, one has
\begin{align*}
A\sum|c_k|^2\leq\Big\|\sum_{i\in[m]}\sum_{k\in \sigma_i} c_k(f_{1k}^i\oplus \cdots \oplus f_{Lk}^i)\Big\|^2 \leq B\sum|c_k|^2
\end{align*}
for all finite scalar sequences $\{c_k\}_{k\in I }$.
\end{defn}

Let us have a look at Example \ref{ex3.7} which shows that if  the family of  frames for atomic spaces  are woven, then the family of vector-valued frames need not be woven. But this is not the case for woven Riesz sequences. More precisely, if a family of Bessel sequences for atomic spaces is a woven Riesz sequence, then  its associated family of vector-valued Bessel sequences is a woven Riesz sequence.
\begin{thm}\label{th4.7}
Suppose $\{F_1^i, \dots, F_L^i, I\}$    is a vector-valued  family for $\mathcal{H}_1\oplus \cdots \oplus\mathcal{H}_L$ such that \break $F_j^i=\{f_{jk}^i\}_{k\in I}$ is a Bessel sequence  for all $i\in [m]$ and $j\in[L]$. If for some $n\in[L]$, the family $\Big\{F_n^i:i\in[m]\Big\}$ is a woven Riesz sequence, then  $\Big\{\{F_1^i, \dots,F_L^i, I\}:i\in[m]\Big\}$ is a woven Riesz sequence.
\end{thm}
\proof
Let $A$ and $B$ be universal Riesz bounds for the family $\Big\{F_n^i:i\in[m]\Big\}$.   Let $\{\sigma_i\}_{i\in[m]}$ be any partition of $I$. Then,  for any  $\{c_k\}_{k\in I}\in\ell^2(I)$, we compute
\begin{align*}
\Big\|\sum_{i\in[m]}\sum_{k\in \sigma_i}c_k\Big(f_{1k}^i\oplus \cdots \oplus f_{Lk}^i\Big)\Big\|^2&=\Big\|\sum_{i\in[m]}\sum_{k\in \sigma_i}c_kf_{1k}^i\Big\|^2+ \cdots +\Big\|\sum_{i\in[m]}\sum_{k\in \sigma_i}c_kf_{Lk}^i\Big\|^2\\
&\geq\Big\|\sum_{i\in[m]}\sum_{k\in \sigma_i}c_kf_{nk}^i\Big\|^2\\
&\geq A\sum_{k\in I}|c_k|^2.
\end{align*}
By Corollary \ref{cor3.8}, $\bigcup\limits_{i\in [m]}\{F_1^i, \dots, F_L^i, \sigma_i\}$ is a Bessel sequence with universal bound $B_o$ (say). Therefore,
\begin{align*}
\Big\|\sum_{i\in[m]}\sum_{k\in \sigma_i}c_k\Big(f_{1k}^i\oplus \cdots \oplus f_{Lk}^i\Big)\Big\|^2\leq B_o\sum_{k\in I}|c_k|^2.
\end{align*}
Hence the family  $\Big\{\{F_1^i, \dots,F_L^i, I\}:i\in[m]\Big\}$ is a woven Riesz sequence with universal bounds $A$ and $B_o$.
\endproof

A vector-valued frame with one component as Riesz basis (for the underlying atomic space) cannot be weave with  vector-valued frames  whose corresponding component is not a Riesz basis.
\begin{prop}
Suppose the  family of vector-valued frames $\Big\{\{F_{1}^i, \dots, F_{L}^i, I\}:i\in[m]\Big\}$ for \break $\mathcal{H}_1\oplus \cdots \oplus\mathcal{H}_L$ is  woven and $F_p^q$ is a Riesz basis for $\mathcal{H}_p$ for some $p\in[L]$ and $q\in[m]$.  Then,  $F_p^i$ is a Riesz basis for $\mathcal{H}_p$ for all $i\in[m]$.   In particular, $\{F_1^i, \dots, F_L^i, I\}$ is a vector-valued Riesz basis for $\mathcal{H}_1\oplus \cdots \oplus\mathcal{H}_L$ for all $i\in[m]$.
\end{prop}
\proof
Since family of vector-valued frames $\Big\{\{F_{1}^i, \dots, F_{L}^i, I\}:i\in[m]\Big\}$ is woven,
 the family \break  $\Big\{F_j^i:i\in[m]\Big\}$ is woven ($j\in[L]$). In particular,  $\Big\{F_p^i:i\in[m]\Big\}$ is woven. By hypothesis,  $F_p^q$ is a Riesz basis for $\mathcal{H}_p$. Therefore,  by Theorem \ref{th2.5},  $F_p^i$ is a Riesz basis for $\mathcal{H}_p$ for all $i\in[m]$.  Hence by Proposition \ref{prop4.4},  the family $\{F_1^i, \dots, F_L^i, I\}$ is a vector-valued Riesz basis for $\mathcal{H}_1\oplus \cdots \oplus\mathcal{H}_L$ for all $i\in[m]$.
\endproof

The following theorem provide sufficient conditions for vector-valued weaving  Riesz sequences in terms of operators on atomic spaces.

\begin{thm}\label{th4.9}
Suppose $\{F_1^i, \dots, F_L^i, \mathbb{N}\}$  $(i\in [m])$  is a vector-valued family for $\mathcal{H}_1\oplus \cdots \oplus\mathcal{H}_L$. Let there exists some $n\in[L]$ such that $F_j^i=\{f_{jk}^i\}_{k\in \mathbb{N}}$ is a Bessel sequence  for $j\in[L]\setminus\{n\}$  $(i\in[m])$. Let $\{\chi_k\}_{k\in \mathbb{N}}$ and $\{e_k\}_{k\in \mathbb{N}}$ be any orthonormal bases for $\mathcal{H}_n$ such that
\begin{align*}
 & f_{nk}^i=\chi_k-\sum\limits_{p \in \mathbb{N}}\alpha_{kp}^iT_p^i(e_k) \ \text{for all} \ i\in[m] \ \text{and for all} \  k\in \mathbb{N },
 \intertext{and}
 &\lambda_i=\sum\limits_{p \in \mathbb{N}}\|T_p^i\|\sup\limits_{k\in \mathbb{N}}|\alpha_{kp}^i|<\frac{1}{3^{\frac{m-1}{2}}} \  \ \text{for all} \ i\in[m],
\end{align*}
where  $T_p^i\in \mathcal{B}(\mathcal{H}_n)$ for each $p \in \mathbb{N}, \ i \in [m]$ and $\alpha_{kp}^i$ are scalars  for all $k, \ p \in \mathbb{N}$, $i \in [m]$ . Then,  $\Big\{\{F_1^i, \dots,F_L^i, I\}:i\in[m]\Big\}$ is a woven Riesz sequence.
\end{thm}
\proof
Let $\{\sigma_i\}_{i\in[m]}$ be any partition of $\mathbb{N}$. Then, for any $N\in\mathbb{N}$ and scalars $c_1, c_2, \dots,c_N$, we compute
\begin{align*}
\Big\|\sum_{i\in[m]}\sum_{k\in[N]\bigcap\sigma_i}c_k(\chi_k-f_{nk}^i)\Big\|
&=\Big\|\sum_{i\in[m]}\sum_{k\in[N]\bigcap\sigma_i}c_k\Big(\sum\limits_{p \in \mathbb{N}}\alpha_{kp}^iT_p^i(e_k)\Big)\Big\|\\
&=\Big\|\sum\limits_{p \in \mathbb{N}}\sum_{i\in[m]}\sum_{k\in[N]\bigcap\sigma_i}c_k\alpha_{kp}^iT_p^i(e_k)\Big\|\\
&\leq\sum\limits_{p \in \mathbb{N}}\Big\|\sum_{i\in[m]}\sum_{k\in[N]\bigcap\sigma_i}c_k\alpha_{kp}^iT_p^i(e_k)\Big\|\\
&\leq\sum\limits_{p \in \mathbb{N}}\sum_{i\in[m]}\Big\|T_p^i\Big(\sum_{k\in[N]\bigcap\sigma_i}c_k\alpha_{kp}^ie_k\Big)\Big\|\\
&\leq\sum\limits_{p \in \mathbb{N}}\sum_{i\in[m]}\|T_p^i\|\Big\|\sum_{k\in[N]\bigcap\sigma_i}c_k\alpha_{kp}^ie_k\Big\|\\
&=\sum\limits_{p \in \mathbb{N}}\sum_{i\in[m]}\|T_p^i\|\Big(\sum_{k\in[N]\bigcap\sigma_i}|c_k\alpha_{kp}^i|^2\Big)^{\frac{1}{2}}\\
&\leq\sum\limits_{p \in \mathbb{N}}\sum_{i\in[m]}\|T_p^i\|\sup\limits_{k\in \mathbb{N}}|\alpha_{kp}^i|\Big(\sum_{k\in[N]\bigcap\sigma_i}|c_k|^2\Big)^{\frac{1}{2}}\\
&=\sum_{i\in[m]}\lambda_i\Big(\sum_{k\in[N]\bigcap\sigma_i}|c_k|^2\Big)^{\frac{1}{2}}\\
&\leq\max\limits_{i\in[m]}\{\lambda_i\}\sum_{i\in[m]}\Big(\sum_{k\in[N]\bigcap\sigma_i}|c_k|^2\Big)^{\frac{1}{2}}\\
&\leq\max\limits_{i\in[m]}\{\lambda_i\}3^{\frac{m-1}{2}}\Big(\sum_{i\in[m]}\sum_{k\in[N]\bigcap\sigma_i}|c_k|^2\Big)^{\frac{1}{2}}\\
&=\lambda\Big(\sum_{k\in[N]}|c_k|^2\Big)^{\frac{1}{2}}, \text{ where } \lambda=3^{\frac{m-1}{2}}\max\limits_{i\in[m]}\{\lambda_i\}<1.
\end{align*}
By Theorem \ref{th2.1}, $\bigcup\limits_{i\in [m]}\{f_{nk}^i\}_{k\in\sigma_i}$ is equivalent to orthonormal basis $\{\chi_k\}_{k\in \mathbb{N}}$. Therefore, there exists a bijective linear operator $U:\mathcal{H}_n\rightarrow \mathcal{H}_n$ such that $U(\chi_k)=f_{nk}^i$ for $k\in\sigma_i$. Thus,  $\bigcup\limits_{i\in [m]}\{f_{nk}^i\}_{k\in\sigma_i}$ is a Riesz basis and hence a frame for $\mathcal{H}_n$. Therefore, the family of frames $\Big\{\{f_{nk}^i\}_{k\in \mathbb{N}}:i\in[m]\Big\}$ is weakly woven and hence woven  (by Remark \ref{rem2.5}). Thus, there exist universal frame bounds (note that the Riesz basis bounds coincide with the frame bounds) for the family  $\Big\{\{f_{nk}^i\}_{k\in \mathbb{N}}:i\in[m]\Big\}$. That is, $\Big\{\{f_{nk}^i\}_{k\in \mathbb{N}}:i\in[m]\Big\}$ is  woven Riesz sequence. Hence  by Theorem \ref{th4.7}, the family $\Big\{\{F_1^i, \dots, F_L^i, \mathbb{N}\}:i\in[m]\Big\}$ is a woven Riesz sequence.
\endproof
To conclude the paper, we  illustrates  Theorem \ref{th4.9} with the following example.
\begin{exa}
Let $L=m=n=2$ and let  $\mathcal{H}_1=\mathcal{H}_2=\ell^2(\mathbb{N})$. Suppose that $\{\chi_k\}_{k\in\mathbb{N}}$ and $\{e_k\}_{k\in\mathbb{N}}$ are  any two orthonormal basis for $\ell^2(\mathbb{N})$.

 Choose $0 < \epsilon<\frac{1}{\sqrt{3}\sum\limits_{k\in\mathbb{N}}1/k^2}$.  Define sequences $\{T_p^i\}_{p \in \mathbb{N}} \subset  \mathcal{B}(\mathcal{H}_2)$ $(i=1, 2)$ by
   \begin{align*}
&T_p^1:\mathcal{H}_2\rightarrow \mathcal{H}_2\text{ by } T_p^1(\{\xi_1, \xi_2, \dots \})=\{0,0,\dots,0,\xi_p, \xi_{p+1}, \dots\},\\
&T_p^2:\mathcal{H}_2\rightarrow \mathcal{H}_2 \text{ by } T_p^2(\{\xi_1, \xi_2, \dots\})= \Big\{\frac{\xi_p}{4},\frac{\xi_{p+1}}{4}, \dots\Big\}.
\end{align*}

Choose $\{\alpha_{kp}^i\}_{k\in \mathbb{N}} \subset \mathbb{K} \ (p \in \mathbb{N}; \ i=1,2)$ as follows:
\begin{align*}
&\{\alpha_{kp}^1\}_{k\in \mathbb{N}}=\Big\{0, \dots, 0,\underbrace{\frac{\epsilon}{p^2}}_{p^{th}-place},\frac{\epsilon}{(p+1)^2}, \dots\Big\}\  (p\in\mathbb{N}),\\
&\{\alpha_{kp}^2\}_{k\in \mathbb{N}}=\Big\{0, \dots,0,\underbrace{\frac{4\epsilon}{(p+1)^2}}_{p+1^{th}-place},\frac{4\epsilon}{(p+2)^2}, \dots\Big\},\  p\in\mathbb{N}.
\end{align*}
Then
\begin{align*}
&\lambda_1=\sum\limits_{p\in \mathbb{N}}\|T_p^1\|\sup\limits_{k\in \mathbb{N}}|\alpha_{kp}^i|\leq\sum\limits_{p\in \mathbb{N}} \frac{\epsilon}{p^2}<\frac{1}{\sqrt{3}},
\intertext{ and }
&\lambda_2=\sum\limits_{p\in \mathbb{N}}\|T_p^2\|\sup\limits_{k\in \mathbb{N}}|\alpha_{kp}^i|\leq\sum\limits_{p\in \mathbb{N}}\frac{1}{4} \frac{4\epsilon}{(p+1)^2}\leq\sum\limits_{p\in \mathbb{N}}\frac{\epsilon}{p^2}<\frac{1}{\sqrt{3}}.
\end{align*}
Define $\{F_1^1,F_2^1,\mathbb{N}\}$ and $\{F_1^2,F_2^2,\mathbb{N}\}$ as follows:
\begin{align*}
&f_{1k}^1 = \Big\{0, \dots,0,\underbrace{1}_{k^{th}-place},1,0,0, \dots\Big\},\  k\in\mathbb{N},\\
&f_{2k}^1=\chi_k-\sum\limits_{p\in \mathbb{N}}\alpha_{kp}^1T_p^1(e_k),\  k\in\mathbb{N},\\
&f_{1k}^2 = \Big\{0, \dots,0,\underbrace{1}_{k^{th}-place},1, 1,0,0, \dots\Big\},\  k\in\mathbb{N},\\
&f_{2k}^2=\chi_k-\sum\limits_{p\in \mathbb{N}}\alpha_{kp}^2T_p^2(e_k),\  k\in\mathbb{N}.
\end{align*}
It is easy to verify that $F_j^i$ is a Bessel sequence for $\mathcal{H}_i$, for all $i\in[m], \ j\in[L]\setminus\{n\}=\{1\}$. Hence,  by  Theorem \ref{th4.9}, the family  $\Big\{\{F_1^i,F_2^i, \mathbb{N}\}:i\in[2]\Big\}$ is a woven Riesz sequence.
\end{exa}

\mbox{}

\end{document}